\documentstyle[twoside]{article}
\title{Evaluations of some terminating hypergeometric ${}_2F_1(2)$ series}
\author{\sc Y. S. Kim$^a$,  A. K. Rathie$^b$ and R. B. Paris$^c$\\
\\
${}^a\!$ Department of Mathematics Education, Wonkwang University, Iksan, Korea\\
E-Mail: yspkim@wonkwang.ac.kr\\
${}^b\!$ Department of Mathematics, Central University of Kerala, Kasaragad 671123,\\
Kerala, India\\
E-Mail: akrathie@cukerala.edu.in\\
${}^c\!$ School of Engineering, Computing and Applied Mathematics,\\
 University of Abertay Dundee, Dundee DD1 1HG, UK\\
E-Mail: r.paris@abertay.ac.uk}
\begin{document}
\def\f#1#2{\mbox{${\textstyle \frac{#1}{#2}}$}}
\def\dfrac#1#2{\displaystyle{\frac{#1}{#2}}}
\def\boldal{\mbox{\boldmath $\alpha$}}
\newcommand{\bee}{\begin{equation}}
\newcommand{\ee}{\end{equation}}
\newcommand{\lam}{\lambda}
\newcommand{\ka}{\kappa}
\newcommand{\al}{\alpha}
\newcommand{\th}{\theta}
\newcommand{\om}{\omega}
\newcommand{\Om}{\Omega}
\newcommand{\fr}{\frac{1}{2}}
\newcommand{\fs}{\f{1}{2}}
\newcommand{\g}{\Gamma}
\newcommand{\br}{\biggr}
\newcommand{\bl}{\biggl}
\newcommand{\ra}{\rightarrow}
\newcommand{\mbint}{\frac{1}{2\pi i}\int_{c-\infty i}^{c+\infty i}}
\newcommand{\mbcint}{\frac{1}{2\pi i}\int_C}
\newcommand{\mboint}{\frac{1}{2\pi i}\int_{-\infty i}^{\infty i}}
\newcommand{\gtwid}{\raisebox{-.8ex}{\mbox{$\stackrel{\textstyle >}{\sim}$}}}
\newcommand{\ltwid}{\raisebox{-.8ex}{\mbox{$\stackrel{\textstyle <}{\sim}$}}}
\renewcommand{\topfraction}{0.9}
\renewcommand{\bottomfraction}{0.9}
\renewcommand{\textfraction}{0.05}
\newcommand{\mcol}{\multicolumn}
\date{}
\maketitle
\begin{abstract}
Explicit expressions for the hypergeometric series ${}_2F_1(-n, a; 2a\pm j;2)$ and ${}_2F_1(-n, a; -2n\\ \pm j;2)$
for positive integer $n$ and arbitrary integer $j$ are obtained with the help of generalizations of Kummer's second and third summation theorems obtained earlier by Rakha and Rathie. Results for $|j|\leq 5$ derived previously using different methods are also obtained as special cases.
\vspace{0.4cm}

 \setcounter{section}{1}
\setcounter{equation}{0}
\vspace{0.3cm}

\end{abstract}
\vspace{0.3cm}

\begin{center}
{\bf 1. \  Introduction}
\end{center}
\setcounter{section}{1}
\setcounter{equation}{0}
\renewcommand{\theequation}{\arabic{section}.\arabic{equation}}
In a problem arising in a model of a biological problem, Samoletov \cite{S} obtained by means of a mathematical induction argument the following sum containing factorials 
\[\sum_{k=0}^n \frac{(-1)^k}{(n-k)!}\,\frac{(2k+1)!!}{k! (k+1)!} = \frac{(-1)^n}{\sqrt{n! (n+1)!}} \left(\sqrt{n+1}\ \frac{(n-1)!!}{n!!}\right)^{(-1)^n},\] 
where throughout $n$ denotes a positive integer and, as usual, 
\[(2n)!!=2\cdot4\cdot 6\cdots=2^n n!,\qquad (2n+1)!!=1\cdot 3\cdot 5\cdots =\frac{(2n+1)!}{2^n n!}.\]
Samoletov also expressed the above sum in the equivalent hypergeometric form
\[{}_2F_1\left[\begin{array}{c}-n, \f{3}{2}\\2\end{array};2\right]=
\left\{\begin{array}{ll}
\frac{\g(\fs n+\fs)}{\surd\pi \g(\fs n+1)} & (n \ \mbox{even})\\
\frac{-\g(\fs n+1)}{\surd\pi \g(\fs n+\f{3}{2})} & (n \ \mbox{odd}) 
\end{array}\right.\]
Subsequently, Srivastava \cite{HMS} pointed out that this result could be easily derived from a known hypergeometric summation formula \cite[Vol.~2, p.~493]{PBM} for ${}_2F_1(-n,a; 2a-1;2)$ with $a=\f{3}{2}$, which is a contiguous result to the well-known summation
\[{}_2F_1\left[\begin{array}{c}-n, a\\2a\end{array};2\right]=\frac{2^n\surd\pi \g(1-a)}{(2a)_n \g(\fs-\fs n) \g(1-a-\fs n)}\qquad (n=0, 1, 2, \ldots).\]

The aim in this note is to obtain explicit expressions for
\bee\label{e11}
{}_2F_1\left[\begin{array}{c}-n, a\\2a\pm j\end{array};2\right]\qquad \mbox{and}\qquad 
{}_2F_1\left[\begin{array}{c}-n, a\\-2n\pm j\end{array};2\right]
\ee
for arbitrary integer $j$. We shall employ the
following generalizations of Kummer's second and
third summation theorems given in \cite{RR} (we correct a misprint in Theorem 6 of this reference). These are respectively
\[{}_2F_1\left[\begin{array}{c} \alpha, \beta\\\fs(\alpha+\beta\pm j+1)\end{array};\frac{1}{2}\right]=
\frac{\surd\pi \g(\fs \alpha+\fs \beta+\fs\pm\fs j)}{\g(\fs \alpha+\fs) \g(\fs \beta+\fs)}\,\frac{\g(\fs \alpha-\fs \beta+\fs\mp\fs j)}{\g(\fs \alpha-\fs \beta+\fs+\fs j)}\]
\bee\label{e12}
\hspace{4cm}\times \sum_{r=0}^j (\mp 1)^r \left(\!\!\begin{array}{c}j\\r\end{array}\!\!\right)\,\frac{(\fs \beta)_{r/2}}{(\fs \alpha+\fs)_{(r-j)/2}}
\ee
and 
\[{}_2F_1\left[\begin{array}{c}\alpha, 1-\alpha\pm j\\\gamma\end{array};\frac{1}{2}\right]=\frac{2^{\pm j}\g(\fs \gamma) \g(\fs\gamma+\fs)}{\g(\fs\gamma+\fs \alpha) \g(\fs\gamma-\fs\alpha+\fs)}\,\frac{\g(\alpha\mp j)}{\g(\alpha+\epsilon_j)}\hspace{3cm}\]
\bee\label{e13}
\hspace{4cm}\times \sum_{r=0}^j(\mp 1)^r \left(\!\!\begin{array}{c}j\\r\end{array}\!\!\right)\,\frac{(\fs\gamma-\fs\alpha)_{r/2}}{(\fs\gamma+\fs\alpha)_{ r/2-\delta_j}} 
\ee
for $j=0, 1, 2, \ldots\,$, where $\epsilon_j=0$ (resp.~$j$), $\delta_j=j$ (resp. 0) when the upper (resp. lower) signs are taken and 
$(a)_k=\g(a+k)/\g(a)$ is the Pochhammer symbol defined for {\it arbitrary} index $k$.
When $j=0$, the summations (\ref{e12}) and (\ref{e13}) reduce to the well-known second and third summation theorems due to Kummer \cite[p.~243]{Sla}
\[{}_2F_1\left[\begin{array}{c} \alpha, \beta\\\fs(\alpha+\beta+1)\end{array};\frac{1}{2}\right]=\frac{\surd\pi \g(\fs \alpha+\fs \beta+\fs)}{\g(\fs \alpha+\fs) \g(\fs \beta+\fs)}\]
and\footnote{In \cite[p.~243]{Sla}, this summation formula is referred to as Bailey's theorem. However, it has been pointed out in \cite{HMS2} that this theorem was originally found by Kummer.}
\[{}_2F_1\left[\begin{array}{c}\alpha, 1-\alpha\\\gamma\end{array};\frac{1}{2}\right]=\frac{\g(\fs \gamma) \g(\fs \gamma+\fs)}{\g(\fs \gamma+\fs \alpha) \g(\fs \gamma-\fs \alpha+\fs)}.\]
In addition, we shall make use of the transformation \cite[(15.8.6)]{DLMF}
\bee\label{e14}
{}_2F_1\left[\begin{array}{c} -n, \beta\\\gamma\end{array};2\right]=\frac{(-2)^n (\beta)_n}{(\gamma)_n}\,{}_2F_1\left[
\begin{array}{c} -n, 1-\gamma-n\\1-\beta-n\end{array};\frac{1}{2}\right].
\ee

Expressions for the series in (\ref{e11}) for arbitrary integer $j$ have recently been obtained by
Chu \cite{C} using a different approach. This involved expressing the series  for $j\neq 0$ as finite sums of ${}_2F_1(2)$ series in (\ref{e11}) with $j=0$.
The cases with $|j|\leq 5$ have also been given by Kim and Rathie \cite{KR} and Kim {\it et al.\/} \cite{KRR}.
An application of the first series in (\ref{e11}) for $j=0, 1, \ldots ,5$ has been discussed in \cite{KRP}.
\vspace{0.6cm}

\begin{center}
{\bf 2. \  Statement of the results}
\end{center}
\setcounter{section}{2}
\setcounter{equation}{0}
\renewcommand{\theequation}{\arabic{section}.\arabic{equation}}
Our principal results are stated in the following two theorems.

\newtheorem{theorem}{Theorem}
\begin{theorem}
Let $n$ be a positive integer, $a$ be a complex parameter and define $j_0=\lfloor \fs j\rfloor$. Then we have
\bee\label{e21a}
{}_2F_1\left[\begin{array}{c} -2n, a\\2a\pm j\end{array};2\right]=
\frac{2^{2n}(\fs)_n}{(2a\pm j)_{2n}} \sum_{r=0}^{j_0} (-1)^r\left(\!\!\begin{array}{c}j\\2r\end{array}\!\!\right) (-n)_r (a+\delta_j)_{n-r}
\ee
and
\bee\label{e21b}
{}_2F_1\left[\begin{array}{c} -2n-1, a\\2a\pm j\end{array};2\right]
=\frac{\pm2^{2n}(\f{3}{2})_n}{(2a\pm j)_{2n+1}} \sum_{r=0}^{j_0} (-1)^r\left(\!\!\begin{array}{c}j\\2r+1\end{array}\!\!\right) (-n)_r (a+\delta_j)_{n-r}
\ee
for $j=0, 1, 2, \ldots\,$, where $\delta_j=j$ (resp.~0) when the upper (resp. lower) signs are taken. 
\end{theorem}

\noindent{\bf Proof}.\ \ From the result (\ref{e14}) we have
\[{}_2F_1\left[\begin{array}{c}-n, a\\2a\pm j\end{array};2\right]=\frac{(-2)^n (a)_n}{(2a\pm j)_n}
\,{}_2F_1\left[\begin{array}{c} -n, 1-2a\mp j-n\\1-a-n\end{array};\frac{1}{2}\right].\]
The hypergeometric function on the right-hand side can be summed by the generalized second Kummer summation theorem (\ref{e12}), where we put $\alpha=1-2a-j-n$ and $\beta=-n$. After some straightforward algebra we obtain
\[{}_2F_1\left[\begin{array}{c}-n, a\\2a+j\end{array};2\right]=\frac{2^n\surd\pi}{(2a+j)_n \g(-\fs n) \g(-\fs n+\fs)}\sum_{r=0}^j(-1)^r\left(\!\!\begin{array}{c}j\\r\end{array}\!\!\right)\frac{\g(-\fs n+\fs r)}{(1-a-j)_{(r-n)/2}}\]
and 
\[{}_2F_1\left[\begin{array}{c}-n, a\\2a-j\end{array};2\right]=\frac{2^n\surd\pi}{(2a-j)_n \g(-\fs n) \g(-\fs n+\fs)}\sum_{r=0}^j\left(\!\!\begin{array}{c}j\\r\end{array}\!\!\right)\frac{\g(-\fs n+\fs r)}{(1-a)_{(r-n)/2}}.\]
Changing $n$ to $2n$ and $2n+1$ and using the properties of the gamma function and
\[(a)_{2n}=2^{2n}(a)_n(a+\fs)_n,\qquad (a)_{2n+1}=2^{2n}  a (\fs a+\fs)_n (\fs a+1)_n,\]
we then find the results stated in the theorem.\ \ \ \ \ $\Box$
\bigskip

We remark that in (\ref{e21b}) the upper limit of summation can be replaced by $\lfloor \fs j\rfloor-1$ when $j$ is even. Also, since $(-n)_r$ vanishes when $r>n$, it is possible to replace the upper summation limit in both (\ref{e21a}) and (\ref{e21b}) by $n$ whenever $n>\lfloor\fs j\rfloor$.

\begin{theorem}
Let $n$ be a positive integer and $a$ be a complex parameter. Then we have
\bee\label{e22a}
{}_2F_1\left[\begin{array}{c}-n, a\\-2n+j\end{array};2\right]=\frac{2^{2n-j} (n-j)!}{(2n-j)!} \sum_{r=0}^j\left(\!\!\begin{array}{c}j\\r\end{array}\!\!\right)(\fs a+\fs -\fs r)_n
\ee
provided $j$ does not lie in the interval $[n+1, 2n]$ where the hypergeometric function on the left-hand side of (\ref{e22a}) is, in general, not defined,
and
\bee\label{e22b}
{}_2F_1\left[\begin{array}{c}-n, a\\-2n-j\end{array};2\right]=\frac{2^{2n+j} n!}{(2n+j)!} \sum_{r=0}^j
(-1)^r \left(\!\!\begin{array}{c}j\\r\end{array}\!\!\right)(\fs a+\fs-\fs r)_{n+j}
\ee
for $j=0, 1, 2, \ldots\,$. When $j\geq 2n+1$ in (\ref{e22a}), the ratio of factorials $(n-j)!/(2n-j)!$ can be replaced by $(-1)^n (j-2n-1)!/(j-n-1)!$.
\end{theorem}

\noindent{\bf Proof}.\ \ From the result (\ref{e14}) we have
\[{}_2F_1\left[\begin{array}{c}-n, a\\-2n\mp j\end{array};2\right]=\frac{2^n (a)_n (n\pm j)!}{(2n\pm j)_n}
\,{}_2F_1\left[\begin{array}{c} -n, 1+n\pm j\\1-a-n\end{array};\frac{1}{2}\right]\]
when the parameters are such that the hypergeometric functions make sense. The hypergeometric function on the right-hand side can be summed by the generalized third Kummer theorem (\ref{e13}), where we put $\alpha=-n$ and $\gamma=1-a-n$. Some straightforward algebra using the properties of the gamma function then yields the results 
(\ref{e22a}) and (\ref{e22b}) in the theorem.\ \ \ \ \ $\Box$
\bigskip

The sums on the right-hand sides of (\ref{e22a}) and (\ref{e22b}) can be written in an alternative form involving
just two Pochhammer symbols containing the index $n$ by making use of the result
\[(\alpha-r)_n=\frac{(\alpha)_n (1-\alpha)_r}{(1-\alpha-n)_r}\]
for positive integers $r$ and $n$. Then we find, with $j_0=\lfloor \fs j\rfloor$,
\[{}_2F_1\left[\begin{array}{c}-n, a\\-2n+j\end{array};2\right]\hspace{10cm}\]
\bee\label{e23a}
\hspace{1cm}=\frac{2^{2n-j} (n-j)!}{(2n-j)!}\left\{(\fs a+\fs)_n\sum_{r=0}^{j_0}\left(\!\!\begin{array}{c}j\\2r\end{array}\!\!\right)A_r(n,0)
+(\fs a)_n\sum_{r=0}^{j_0}\left(\!\!\begin{array}{c}j\\2r+1\end{array}\!\!\right)B_r(n,0)\right\}
\ee
and
\[{}_2F_1\left[\begin{array}{c}-n, a\\-2n-j\end{array};2\right]\hspace{10cm}\]
\bee\label{e23b}
\hspace{0.8cm}=\frac{2^{2n+j}n!}{(2n+j)!}\left\{(\fs a+\fs )_{n+j}\sum_{r=0}^{j_0}\left(\!\!\begin{array}{c}j\\2r\end{array}\!\!\right)A_r(n,j)-(\fs a)_{n+j}
\sum_{r=0}^{j_0}\left(\!\!\begin{array}{c}j\\2r+1\end{array}\!\!\right)B_r(n,j)\right\},
\ee
where
\[A_r(n,j):=\frac{(\fs-\fs a)_r}{(\fs-\fs a-n-j)_r},\qquad B_r(n,j):=\frac{(1-\fs a)_r}{(1-\fs a-n-j)_r}.\]
Again, when $j$ is even, the upper summation limit in the second sums in (\ref{e23a}) and (\ref{e23b}) can be replaced by $j_0-1$, if so desired.
\vspace{0.6cm}

\begin{center}
{\bf 3. \  Special cases}
\end{center}
\setcounter{section}{3}
\setcounter{equation}{0}
\renewcommand{\theequation}{\arabic{section}.\arabic{equation}}
If we set $0\leq j\leq 5$ in (\ref{e21a}) we obtain the following summations:
\bee{}_2F_1\left[\begin{array}{c}-2n, a\\2a\end{array};2\right]=\frac{(\fs)_n}{(a+\fs)_n}={}_2F_1\left[\begin{array}{c}-2n, a\\2a+1\end{array};2\right],\ee
\bee{}_2F_1\left[\begin{array}{c}-2n, a\\2a+2\end{array};2\right]=\frac{(\fs)_n}{(a+\f{3}{2})_n}\left(1+\frac{2n}{a+1}\right),\ee
\bee{}_2F_1\left[\begin{array}{c}-2n, a\\2a+3\end{array};2\right]=\frac{(\fs)_n}{(a+\f{3}{2})_n}\left(1+\frac{4n}{a+2}\right),\ee
\bee{}_2F_1\left[\begin{array}{c}-2n, a\\2a+4\end{array};2\right]=\frac{(\fs)_n}{(a+\f{5}{2})_n}\left(1+\frac{8n}{a+2}+\frac{16n(n-1)}{(a+2)(a+3)}\right),\ee
\bee{}_2F_1\left[\begin{array}{c}-2n, a\\2a+5\end{array};2\right]=\frac{(\fs)_n}{(a+\f{5}{2})_n}\left(1+\frac{12n}{a+3}+\frac{16n(n-1)}{(a+3)(a+4)}\right)\ee
and
\bee{}_2F_1\left[\begin{array}{c}-2n, a\\2a-1\end{array};2\right]=\frac{(\fs)_n}{(a-\fs)_n},\ee
\bee{}_2F_1\left[\begin{array}{c}-2n, a\\2a-2\end{array};2\right]=\frac{(\fs)_n}{(a-\fs)_n}\left(1+\frac{2n}{a-1}\right),\ee
\bee{}_2F_1\left[\begin{array}{c}-2n, a\\2a-3\end{array};2\right]=\frac{(\fs)_n}{(a-\f{3}{2})_n}\left(1+\frac{4n}{a-1}\right),\ee
\bee{}_2F_1\left[\begin{array}{c}-2n, a\\2a-4\end{array};2\right]=\frac{(\fs)_n}{(a-\f{3}{2})_n}\left(1+\frac{8n}{a-2}+\frac{8n(n-1)}{(a-1)(a-2)}\right),\ee
\bee{}_2F_1\left[\begin{array}{c}-2n, a\\2a-5\end{array};2\right]=\frac{(\fs)_n}{(a-\f{5}{2})_n}\left(1+\frac{12n}{a-2}+\frac{16n(n-1)}{(a-1)(a-2)}\right).\ee

Similarly, if we set $0\leq j\leq 5$ in (\ref{e21b}) we obtain the following summations:
\bee{}_2F_1\left[\begin{array}{c}-2n-1, a\\2a\end{array};2\right]=0,\ee
\bee{}_2F_1\left[\begin{array}{c}-2n-1, a\\2a+1\end{array};2\right]=\frac{(\f{3}{2})_n}{(2a+1)(a+\f{3}{2})_n},\ee
\bee{}_2F_1\left[\begin{array}{c}-2n-1, a\\2a+2\end{array};2\right]=\frac{2(\f{3}{2})_n}{(2a+2)(a+\f{3}{2})_n},\ee
\bee{}_2F_1\left[\begin{array}{c}-2n-1, a\\2a+3\end{array};2\right]=\frac{(\f{3}{2})_n}{(2a+3)(a+\f{5}{2})_n}\left(3+\frac{4n}{a+2}\right),\ee
\bee{}_2F_1\left[\begin{array}{c}-2n-1, a\\2a+4\end{array};2\right]=\frac{(\f{3}{2})_n}{(2a+4)(a+\f{5}{2})_n}\left(4+\frac{8n}{a+3}\right),\ee
\bee{}_2F_1\left[\begin{array}{c}-2n-1, a\\2a+5\end{array};2\right]=\frac{(\f{3}{2})_n}{(2a+5)(a+\f{7}{2})_n}\left(5+\frac{20n}{a+3}+\frac{16n(n-1)}{(a+3)(a+4)}\right)\ee
and
\bee{}_2F_1\left[\begin{array}{c}-2n-1, a\\2a-1\end{array};2\right]=-\frac{(\f{3}{2})_n}{(2a-1)(a+\fs)_n},\ee
\bee{}_2F_1\left[\begin{array}{c}-2n-1, a\\2a-2\end{array};2\right]=-\frac{2(\f{3}{2})_n}{(2a-2)(a-\fs)_n},\ee
\bee{}_2F_1\left[\begin{array}{c}-2n-1, a\\2a-3\end{array};2\right]=-\frac{(\f{3}{2})_n}{(2a-3)(a-\fs)_n}\left(3+\frac{4n}{a-1}\right),\ee
\bee{}_2F_1\left[\begin{array}{c}-2n-1, a\\2a-4\end{array};2\right]=-\frac{(\f{3}{2})_n}{(2a-4)(a-\f{3}{2})_n}\left(4+\frac{8n}{a-1}\right),\ee
\bee{}_2F_1\left[\begin{array}{c}-2n-1, a\\2a-5\end{array};2\right]=
-\frac{(\f{3}{2})_n}{(2a-5)(a-\f{3}{2})_n}\left(5+\frac{20n}{a-2}+\frac{16n(n-1)}{(a-1)(a-2)}\right).\ee

Finally, from (\ref{e23a}) and (\ref{e23b}) we obtain:
\bee{}_2F_1\left[\begin{array}{c}-n, a\\-2n\end{array};2\right]=\frac{2^{2n}n!}{(2n)!} (\fs a+\fs)_n=\frac{(\fs a+\fs )_n}{(\fs)_n},\ee
\bee{}_2F_1\left[\begin{array}{c}-n, a\\-2n+1\end{array};2\right]=\frac{2^{2n-1}(n-1)!}{(2n-1)!}\{(\fs a+\fs)_n+(\fs a)_n\},\ee
\bee{}_2F_1\left[\begin{array}{c}-n, a\\-2n-1\end{array};2\right]= \frac{2^{2n+1}n!}{(2n+1)!}\,(\fs a+\fs)_{n+1},\ee
\bee{}_2F_1\left[\begin{array}{c}-n,a\\-2n+2\end{array};2\right]=\frac{2^{2n-1}(n-2)!}{(2n-2)!}\left\{\frac{1-a-n}{1-a-2n}\,(\fs a+\fs)_n+(\fs a)_n\right\},\ee
\bee{}_2F_1\left[\begin{array}{c}-n, a\\-2n-2\end{array};2\right]=\frac{2^{2n+3}n!}{(2n+2)!}\left\{
\frac{(1-a-n-j)}{1-a-2n-2j}\,(\fs a+\fs)_{n+2}-(\fs a)_{n+2}\right\}\ee
and so on. 

The above evaluations agree with those given in \cite{C, KR}, although presented in a different format;
the results (3.1) and (3.6), together with (3.11), (3.12) and (3.17), are also recorded in \cite{Prud} in another form.
\vspace{0.6cm}

\begin{center}
{\bf 4. \  An application of Theorem 1}
\end{center}
\setcounter{section}{4}
\setcounter{equation}{0}
\renewcommand{\theequation}{\arabic{section}.\arabic{equation}}
Kummer's second theorem applied to the confluent hypergeometric function ${}_1F_1$ is \cite[p.~12]{Sla2}
\bee\label{e41}
e^{-x/2} {}_1F_1\left[\begin{array}{c}a\\2a\end{array};x\right]={}_0F_1\left[\begin{array}{c}-\\a+\fs\end{array};\frac{x^2}{16}\right]=(\f{1}{4}x)^{\frac{1}{2}-a} \g(a+\fs) I_{a-\frac{1}{2}}(\fs x),
\ee 
where $I_\nu(z)$ denotes modified Bessel function of the first kind.
We now show how the result in Theorem 1 can be used to derive a generalization of the above theorem for the functions
\[e^{-x/2}{}_1F_1\left[\begin{array}{c}a\\2a\pm j\end{array};x\right]\]
for arbitrary integer $j$.

We have upon series expansion
\[e^{-x/2}{}_1F_1\left[\begin{array}{c}a\\2a\pm j\end{array};x\right]=\sum_{n=0}^\infty \frac{(-\fs x)^n}{n!} 
\sum_{m=0}^\infty \frac{(a)_m}{(2a\pm j)_m}\,\frac{x^m}{m!}
=\sum_{n=0}^\infty\sum_{m=0}^\infty\frac{(-1)^n(a)_m x^{m+n}}{2^n (2a\pm j)_m m! n!}.
\]
Making the change of summation index $n\ra n-m$ and using the fact that $(n-m)!=(-1)^m m!/(-n)_m$, we find
\[e^{-x/2}{}_1F_1\left[\begin{array}{c}a\\2a\pm j\end{array};x\right]
=\sum_{n=0}^\infty \sum_{m=0}^n\frac{(-1)^n(a)_m (-n)_mx^n}{2^{n-m}(2a\pm j)_m m! n!}
\]
\[=\sum_{n=0}^\infty\frac{(-\fs x)^n}{n!}\,{}_2F_1\left[\begin{array}{c}-n, a\\2a\pm j\end{array};2\right].\]

Separation of the above sum into even and odd $n$, use of the evaluation of the ${}_2F_1(2)$ series given in (\ref{e21a}) and (\ref{e21b}) followed by inversion of the order of summation then leads to the result
\[e^{-x/2}{}_1F_1\left[\begin{array}{c}a\\2a\pm j\end{array};x\right]
=\sum_{r=0}^{j_0}(-1)^r\left(\!\!\begin{array}{c}j\\2r\end{array}\!\!\right) \sum_{n=0}^\infty \frac{(-n)_r (a+\delta_j)_{n-r} x^{2n}}{2^{2n} (2a\pm j)_{2n} n!}\hspace{3cm}\]
\bee\label{e42}
\hspace{4cm}\mp \sum_{r=0}^{j_0}(-1)^r\left(\!\!\begin{array}{c}j\\2r+1\end{array}\!\!\right)
\sum_{n=0}^\infty\frac{(-n)_r (a+\delta_j)_{n-r} x^{2n+1}}{2^{2n+1} (2a\pm j)_{2n+1} n!},
\ee
where $j_0$ and $\delta_j$ are defined in Theorem 1. When $j=0$ it is easily seen that (\ref{e42}) reduces to (\ref{e41}).

The result (\ref{e42}) is given in a different form in terms of modified Bessel functions in \cite[Vol.~3, p.~579]{PBM}.

\vspace{0.8cm}

\noindent{\bf Acknowledgement:}\ \ \ One of the authors (YSK) acknowledges the support of the Wonkwang University Research Fund (2014).

\end{document}